\documentclass[12pt]{article}
\usepackage{bm}
\usepackage{float}
\usepackage[spanish,activeacute]{babel}
\usepackage{graphicx, mathrsfs,amssymb,amsmath,amsfonts,amsbsy,amsthm,latexsym}
\begin{document}

\newcommand{\Dh}{\hbox{\bf D}}
\newcommand{\Xh}{\hbox{\bf X}}
\newcommand{\Fh}{\hbox{\bf F}}
\newcommand{\nablaslh}{\nabla\raise2pt\hbox{\kern-9pt\slash\kern3pt}}

\title{\vskip.5cm
Quadrature formulas for the Laplace and Mellin transforms\\
\vskip.5cm}
\author{Rafael G. Campos and Francisco Mej\'{\i}a\\
Facultad de Ciencias F\'{\i}sico-Matem\'aticas,\\
Universidad Michoacana, \\
58060, Morelia, Mich., M\'exico.\\
\hbox{\small rcampos@umich.mx, fdiaz@fismat.umich.mx}\\
}
\date{}
\maketitle
{\vskip.3cm
\noindent MSC: 44A10, 65D32, 33C45
}\\
\vspace*{1.5truecm}
\begin{center} Abstract \end{center}
A discrete Laplace transform and its inversion formula are obtained by using a quadrature of the continuous Fourier transform which is given in terms of Hermite polynomials and its zeros. This approach yields a convergent discrete formula for the two-sided Laplace transform if the function to be transformed falls off rapidly to zero and satisfy certain conditions of integrability, achieving convergence also for singular functions. The inversion formula becomes a quadrature formula for the Bromwich integral. This procedure also yields a quadrature formula for the Mellin transform and its corresponding inversion formula that can be generalized straightforwardly for functions of several variables.
\vskip1.5cm
\section{Introduction}\label{intro}
It is commonly accepted that the problem of obtaining a discrete formula for the Laplace transform
\begin{equation}\label{tl}
g(s)=\int_0^\infty e^{-st}f(t)dt
\end{equation}
of a function $f(t)$ it is not so complicated as the inverse problem. This is due to the fact that the problem concerning the inversion of the Laplace transform is an ill-posed problem \cite{Cop90}. Many papers have been written on this subject (see \cite{Val01}, which contains a list of references) but the techniques used to obtain an inversion formula can be classified \cite{Aba04} into only four main groups: those that use Fourier series, Laguerre functions, Gaver functionals, and the ones that discretize the Bromwich contour. The inversion formula for the two-sided Laplace transform presented here belongs to the last group and it is based on a quadrature of the integral Fourier transform \cite{Cam92, Cam95}. This quadrature formula is given in terms of a matrix $N\times N$ whose elements are constructed from the $N$ zeros of the Hermite polynomial $H_N(t)$ and has order $O(1/N)$ if the function to be transformed is square-integrable in $(-\infty,\infty)$ and satisfy certain conditions of integrability \cite{Cam95}. The aim of this paper is to show that a simple and straightforward adaptation of such a formula yields a discrete two-sided Laplace transform with an easy-to-compute inversion formula corresponding to a quadrature of the Bromwich integral, and a discrete Mellin transform and its inversion formula. All of these discrete transforms can be generalized easily to the case of several variables.
\\
\section{A Discrete Laplace transform}
Firstly, we reformulate the procedure followed in \cite{Cam92, Cam95} to obtain a quadrature formula for the integral Fourier transform  yielding a discrete Fourier transform. Proofs and further applications can be found in these references.\\
Let us consider the set of functions $u_n(t)=\exp(-t^2/2)H_n(t)$, $n=0,1,\ldots$, where $H_n(t)$ is the $n$th Hermite polynomial. This set is closed in  $L^2(-\infty,\infty)$ \cite{Sze75} and their elements are related by the recurrence equation $u_{n+1}(t)+2nu_{n-1}(t)=2tu_n(t)$, which can be written as the eigenvalue problem $\mathbb T\,\mathbb U= t\mathbb U$, $-\infty< t< \infty$, where ${\mathbb T}_{nk}=\delta_{n+1,k}/2+(n-1)\delta_{n,k+1}$, $n,k=1,2,\ldots$,  and $\mathbb U$ 
is the vector whose $n$th entry is $u_{n-1}(t)$.\\
The Fourier transform of $u_n(t)$, denoted by  $v_n(\omega)$, is given by
\begin{equation}\label{tfu}
v_n(\omega)=\int_{-\infty}^{\infty}e^{-i \omega t} u_n(t)dt=\sqrt {2\pi}(-i)^nu_n(\omega),
\end{equation}
and satisfy the recurrence equation $v_{n+1}(\omega)-2nv_{n-1}(\omega)=-2i \omega\,v_n(\omega)$, which can be written in the matrix form ${\mathbb W}\,\mathbb V= -i \omega \mathbb V$, $-\infty< \omega< \infty$, where ${\mathbb W}_{nk}=\delta_{n+1,k}/2-(n-1)\delta_{n,k+1}$, $n,k=1,2,\ldots$, and $\mathbb V$ is the vector whose $n$th entry is $v_{n-1}(\omega)$. By writing the recurrence equations as matrix equations we can consider the eigenproblems corresponding to the principal submatrices of order $N$ of ${\mathbb T}$ and ${\mathbb W}$ to generate sequences of $N$-dimensional vectors $U$ y $V$ converging to $\mathbb U$ and $\mathbb V$ respectively when $N\to \infty$ and in this way, to generate approximations to the functions $u_n(t)$ and $v_n(\omega)$. First let us note that the diagonal matrix ${\mathbb S}$ whose elements are given by ${\mathbb S}_{jk}=\sqrt{2^{k-1} (k-1)!}\delta_{jk}$, generates a symmetric matrix ${\mathbb S}^{-1}{\mathbb T}{\mathbb S}$ and a skew-symmetric matrix ${\mathbb S}^{-1}{\mathbb W}{\mathbb S}$ whose principal submatrices of order $N$, denoted by $T$ and $W$, have elements given by $T_{nk}=\sqrt{n/2} \delta_{n+1,k}+\sqrt{(n-1)/2} \delta_{n,k+1}$ and $W_{nk}=\sqrt{n/2} \delta_{n+1,k}-\sqrt{(n-1)/2} \delta_{n,k+1}$, respectively. Now let us consider the finite eigenproblems of $T$ and $W$:
\[
TU_k=t_k U_k,\quad WV_k=\omega_k V_k,\qquad k=1,2,\ldots, N.
\]
The above recurrence equations and the Christoffel-Darboux formula can be used to find the functional form of the eigenvectors, and also to show that the eigenvalues $t_k$ and $\omega_k$ are both zeros of $H_N(x)$. Thus, the $n$th entries of  the eigenvectors $U_k$ and $V_k$ are given by
\begin{equation}\label{ecUVnk}
U_{nk}=\varphi_{n-1}(t_k),\quad V_{nk}=(-i)^{n-1}\varphi_{n-1}(\omega_k),\qquad n=1,\ldots,N,
\end{equation}
where
\[
\varphi_m(x)=\sqrt{     \frac{ (N-1)!2^{N-m-1}}{Nm!}   }\frac{H_m(x)}{H_{N-1}(x)}.
\]
By construction, $T$ and $W$ approach ${\mathbb T}$ and ${\mathbb W}$ respectively  when  $N\to\infty$. Therefore, in this limit, the $n$th elements of $U_k$ and $V_k$ approach $u_n(t_k)$ and $v_n(\omega_k)$ respectively, up to a constant factor. Since $v_n(\omega)$  is the Fourier transform of $u_n(t)$, the linear transformation $F$ which yields the vector $V_{nk}$, $k=1,\ldots,N$, when it is applied to $U_{nk}$, $k=1,\ldots,N$, corresponds to a discretization of the Fourier transform. This transformation is determined by the matrices $U$ and $V$, whose $k$th columns are just $U_k$ and $V_k$ respectively. Since $F$ satisfies the relation $V^T=FU^T$ between the transpose matrices $V^T$ and $U^T$, we get
\begin{equation}\label{ecdft}
F=V^TU.
\end{equation}
The elements of the unitary and symmetric matrix $F$ 
\begin{equation}\label{Fjk}
F_{kj}=\frac{2^{N-1}(N-1)!}{N{H_{N-1}(t_j)H_{N-1}(\omega_k)}} \sum_{n=0}^{N-1}\frac{(-i)^n}{2^nn!}H_n(t_j)H_n(\omega_k),
\end{equation}
satisfy 
\[
F_{kj}=\frac{\Delta t}{\sqrt{2 \pi}}(-1)^{j+k}e^{-i t_j\omega_k}+{\mathcal O}(1/N),
\]
for bounded $t_j$ and $\omega_k$. Here, $\Delta t=t_{j+1}-t_j=\pi/\sqrt{2N}$ is the Riemann measure that yields the quadrature formula
\begin{equation}\label{cuadtf}
\int_{-\infty}^{\infty}e^{-i \omega_k t} f(t)dt=\int_{-\infty}^{\infty}e^{-s_k t} f(t)dt=\sqrt {2\pi}\sum_{j=1}^N (-1)^{j+k} F_{kj}f(t_j)+{\mathcal O}(1/N)
\end{equation}
for the integral Fourier transform of $f(t)$ evaluated at $\omega_k$ and for the two-sided Laplace transform of $f(t)$ evaluated at $s_k=i\omega_k$. The order of this formula holds whenever $f(t)$ satisfies certain conditions of integrability \cite{Cam95}. If furthermore$f(t)$ is a causal function
\[
f(t)=\begin{cases} h(t),& t\ge 0\\ 0,& t<0,\end{cases}
\]
equation (\ref{cuadtf}) becomes a discrete formula for the Laplace transform of $h(t)$
\begin{equation}\label{cuadlap}
g(s_k)=\int_0^{\infty}e^{-s_k t} h(t)dt=\sum_{j=1}^N L_{kj}f(t_j)+{\mathcal O}(1/N),
\end{equation}
where $s_k=i \omega_k$ and 
\begin{equation}\label{lkjlap}
L_{kj}=\sqrt {2\pi} (-1)^{j+k}F_{kj}. 
\end{equation}
The generalization of this discrete transform to several variables is straightforward. Let $g(s^1,s^2,\ldots,s^n)$ be the $n$-dimensional two-sided Laplace transform of $f(t^1,t^2,\ldots,t^n)$, i.e.,
\[
g(s^1,s^2,\ldots,s^n)=\int_{-\infty}^{\infty}e^{-s\cdot t} f(t^1,t^2,\ldots,t^n)dt^1 dt^2\cdots dt^n,
\]
where $s=(s^1,s^2,\ldots,s^n)$ and $t=(t^1,t^2,\ldots,t^n)$. Then, the corresponding discrete transform if given by the matrix
\begin{equation} \label{seis}
{\mathbf L}=L_n\otimes\cdots\otimes L_l\otimes\cdots\otimes L_1
\end{equation}
in which the entries of $L_l$ are built out of $N_l$ Hermite zeros lying on the $l$th direction and the approximant $\tilde{\mathbf g}$ to $g(s^1,s^2,\ldots,s^n)$ is obtained through the product 
\begin{equation}\label{dltmd}
\tilde{\mathbf g}=\mathbf L{\mathbf f},
\end{equation}
where $\mathbf L$ is the matrix defined in (\ref{seis}), $\mathbf f$ is the vector whose components are given and ordered by
\begin{equation}\label{orden}
f_r=f(t^1_{j_1},t^2_{j_2},\cdots,t^n_{j_n}).
\end{equation}
The index $r$ is related to the others by $r=j_1+(j_2-1)N_1+(j_3-1)N_1N_2+\cdots+(j_n-1)\prod_{l=1}^{n-1}N_l$, where $j_l=1,2,\ldots,N_l$. The component $\tilde{g}\,_r$ of the vector $\mathbf g$ is the approximation to the exact transform $g(s^1_{j_1},s^2_{j_2},\cdots,s^n_{j_n})$ where $s^l_{j_l}=i \omega^l_{j_l}$, $\omega^l_{j_l}=t^l_{j_l}$ y $l=1,2\ldots, n$.\\
Notice that this approach on the discrete Laplace transform put the direct and inverse problem on the same footing since $F^{-1}=F^\dagger$. Thus, the elements of $L^{-1}$ can be computed directly by
\[
L^{-1}_{jk}=(-1)^{j+k} F_{jk}^*/\sqrt {2\pi},
\] 
where $^*$ means complex conjugation. By applying $L^{-1}$ to the vector whose elements are the values of $g(s)$ evaluated at $s_k=i \omega_k$, we get an approximation to the values of $f(t)$ at $t_j$. In other words, we obtain a discretization of the Bromwich integral 
\begin{equation}\label{finvlap}
f(t_j)=\frac{1}{2\pi i}\int_{-i\infty}^{i\infty}e^{t_j\,s} g(s)ds=\sum_{j=1}^N L^{-1}_{jk}g(s_k)+{\mathcal O}(1/N),\quad s_k=i \omega_k,
\end{equation}
in which the contour of integration is the imaginary axis and the singularities of $g(s)$ lie on the left of this line. The real part and the imaginary part of $g(s)$ should satisfy the conditions on integrability given in \cite{Cam95}. The extension to several variables is obvious. Now we have
\begin{equation}\label{invdltmd}
\tilde{\mathbf f}=\mathbf L^{-1}{\mathbf g},
\end{equation}
where
\begin{equation}\label{tdlimd}
\mathbf L^{-1}=L_n^{-1}\otimes\cdots\otimes L_l^{-1}\otimes\cdots\otimes L_1^{-1},
\end{equation}
$\mathbf g$ is the function $g(s^1,s^2,\cdots,s^n)$ evaluated at  $s^l_{j_l}=i \omega^l_{j_l}$, $l=1,2\ldots, n$, and $\tilde{\mathbf f}$ is the approximant to $f(t^1,t^2,\ldots,t^n)$.
\\
\section{A discrete Mellin transform}
Since the Mellin transform
\[
g_M(s)=\int_0^\infty x^{s-1}f(x)dx
\]
is a two-sided Laplace transform under the transformation $x=\exp(-t)$, the discrete Laplace transform $\mathbf L$ defined in (\ref{seis}) yields a discretization of the multidimensional Mellin transform evaluated on the imaginary axis of each variable $s^l$, $l=1,\ldots,n$. Thus, we have that if $f(x^1,x^2,\ldots,x^n)$ is the function to be transformed, the pair of discrete multidimensional Mellin transforms are given by the formulas
\begin{equation}\label{tdsm}
\tilde{\mathbf g}_M=\mathbf L{\mathbf f}_t,\qquad \tilde{\mathbf f}_t=\mathbf L^{-1}{\mathbf g}_M,
\end{equation}
in which ${\mathbf f}_t$ is the vector whose elements are given by
\[
f_t(t^1,t^2,\ldots,t^n)=f(\exp(-t^1),\exp(-t^2),\ldots,\exp(-t^n))
\] 
and ordered according to (\ref{orden}).  It should be noticed that in the inverse formula, the vector $\tilde{\mathbf f}_t$ approaches $f_t(t^k)$ instead of $f(x^k)$.
\section{Examples}
In this section we perform some numerical calculations to show the accuracy of the above discrete Laplace and Mellin transforms. We present two singular  cases (the first and third examples) for which the discrete transforms yield convergent results. For such cases the necessary conditions to get the order $O(1/N)$ are not fulfilled, therefore, the order of convergence is estimated numerically in the next section.
\subsection{Discrete Laplace transforms}
As a first example, we compute the numerical inversion of 
\begin{equation}\label{ejpri}
g(s)=2\sum_{k=1}^n \cosh(ks),
\end{equation}
which is the two-sided Laplace transform of a train of $2n$ delta functions centered at the integers $\pm 1,\pm 2,\ldots,\pm n$. This problem resembles the numerical inversion of the partition function of the quantum harmonic oscillator, a typical test problem.\\
In order to approximate the inverse transform of (\ref{ejpri}), the number $N$ of zeros of $H_N(t)$ should be greater than $n^2/2$ because in this way the interval $[-n,n]$ is contained in $(-\sqrt{2N+1},\sqrt{2N+1})$, which is the interval where the Hermite zeros lie. The application of $L^{-1}$ to the vector $g$ whose elements are the values of (\ref{ejpri}) at the Hermite zeros on the imaginary axis yields the interpolated set of points shown in Figure 1. The result is a function showing the typical features of a sum of delta functions centered at integer values.
\vskip1cm
\hbox to \textwidth{\hfill\scalebox{0.76}{\includegraphics{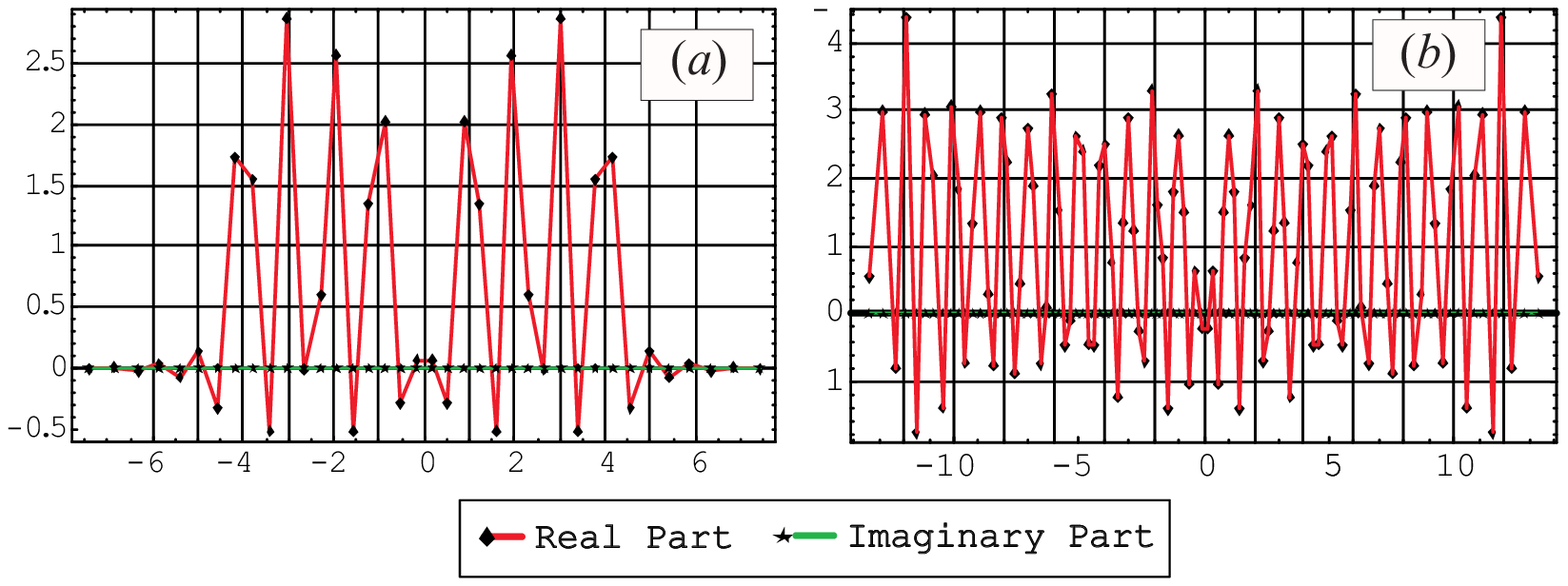}}\hfill}
\begin{center}
\begin{minipage}{13cm}
{\small 
Figure 1: Two-sided numerical inversion of the Laplace transformed function (\ref{ejpri}). In ({\it a}) 40 Hermite zeros have been used for $n=4$. In ({\it b}) 100 Hermite zeros have been used for $n=13$. The maxima of the real parts are centered at the corresponding integers and the imaginary parts are zero.}
\end{minipage}
\end{center}
As a second example, we take the function $h(t)=\exp(-t) \sin(t)$, $0<t<\infty$, whose one-sided Laplace transform is given by $g(s)=1/[(s+1)^2+1]$. According to 
(\ref{cuadlap}), $h(t)$ should be substituted by the causal function 
\begin{equation}\label{ejlapd}
f(t)=\begin{cases} \exp(-t) \sin(t),& t\ge 0\\ 0,& t<0,\end{cases}
\end{equation}
in order to obtain the approximated Laplace transform. The application of (\ref{cuadlap}) and (\ref{finvlap}) to the vectors $f$ and $g$ respectively, yields the results displayed in Figure 2. For $N=40$, the relative errors are given by
\[
\frac{\Vert g-\tilde{g}  \Vert_2}{\Vert g \Vert_2}=0.023758,\qquad \frac{\Vert f-\tilde{f}  \Vert_2}{\Vert f \Vert_2}=0.0236836.
\]
It should be reminded that $g$, $\tilde{g}$ and $\tilde{f}$ are complex vectors.
\vskip1cm
\hbox to \textwidth{\hfill\scalebox{0.8}{\includegraphics{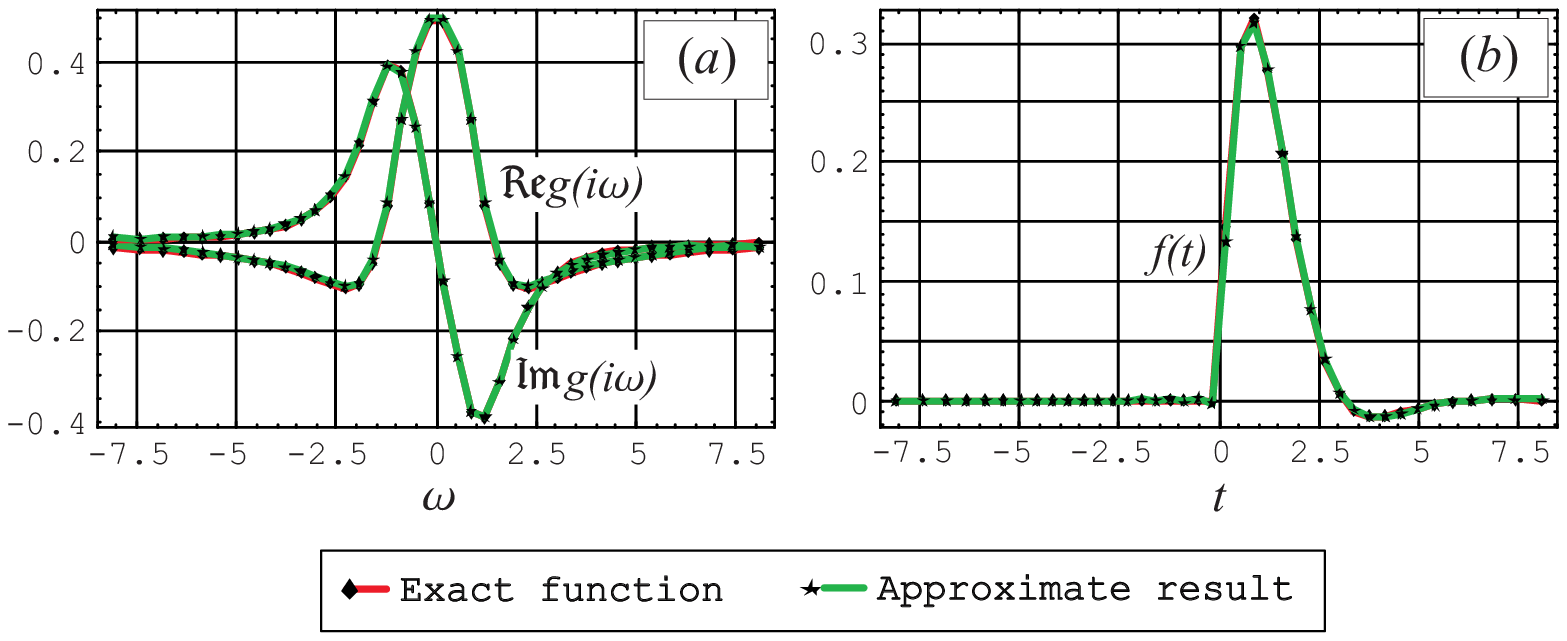}}\hfill}
\begin{center}
\begin{minipage}{13cm}
{\small 
Figure 2: ({\it a}) Exact and approximate Laplace transform of (\ref{ejlapd}). ({\it b}) Inverse transform obtained through (\ref{invdltmd}).
In both cases, 40 Hermite zeros on the imaginary axis were used.}
\end{minipage}
\end{center}
The next examples concern the performance of the discrete Mellin transform. As a first case, we take the singular problem defined by the Mellin transform of 
\begin{equation}\label{ejtmelu}
f(x)=\frac{\sqrt{x}}{1-x},\quad 0<x<\infty.
\end{equation}
The Cauchy principal value of this integral is $-\pi \tan(\pi s)$ and it is displayed in Figure 3, together with the discrete Mellin transforms (\ref{tdsm}). Figure 3($b$) shows the plot of $f(\exp(-t))$ against $t$ instead of $f(x)$ against $x$. The corresponding relative errors are
\[
\frac{\Vert g_M-\tilde{g}_M  \Vert_2}{\Vert g_M \Vert_2}=0.156919,\qquad \frac{\Vert f_t-\tilde{f}_t  \Vert_2}{\Vert f_t \Vert_2}=0.0739943.
\]
\vskip1cm
\hbox to \textwidth{\hfill\scalebox{0.8}{\includegraphics{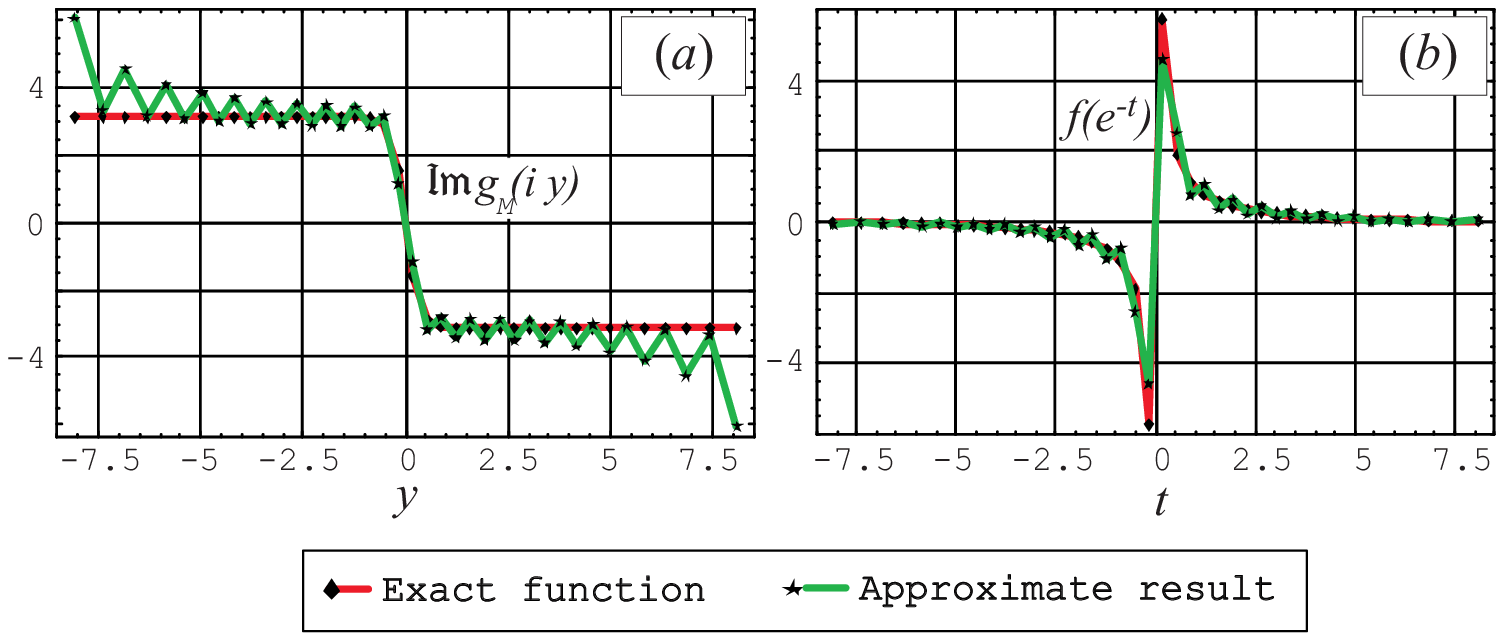}}\hfill}
\begin{center}
\begin{minipage}{13cm}
{\small 
Figure 3: ({\it a}) Exact and approximate Mellin transform of (\ref{ejtmelu}) on the imaginary axis (the real part is zero). ({\it b}) Inverse transform obtained by (\ref{tdsm}). In this case the imaginary part is zero. In both cases 40 Hermite zeros were used.}
\end{minipage}
\end{center}
As a final example, we consider the function
\begin{equation}\label{ejcmel}
f(x)=\exp(-\frac{x}{\sqrt{2}}) \sin(\frac{x}{\sqrt{2}}), \quad  0<x<\infty,
\end{equation}
whose Mellin transform is $\sin(\pi s/4)\Gamma(s)$. Figure 4 shows the output of the discrete transforms (\ref{tdsm}). Again, $f(\exp(-t))$ is plotted against $t$ in Figure 4($b$). The relative errors are
\[
\frac{\Vert g-\tilde{g}  \Vert_2}{\Vert g \Vert_2}=0.00702041,\qquad \frac{\Vert f-\tilde{f}  \Vert_2}{\Vert f \Vert_2}=0.00701767.
\]
\vskip1cm
\hbox to \textwidth{\hfill\scalebox{0.8}{\includegraphics{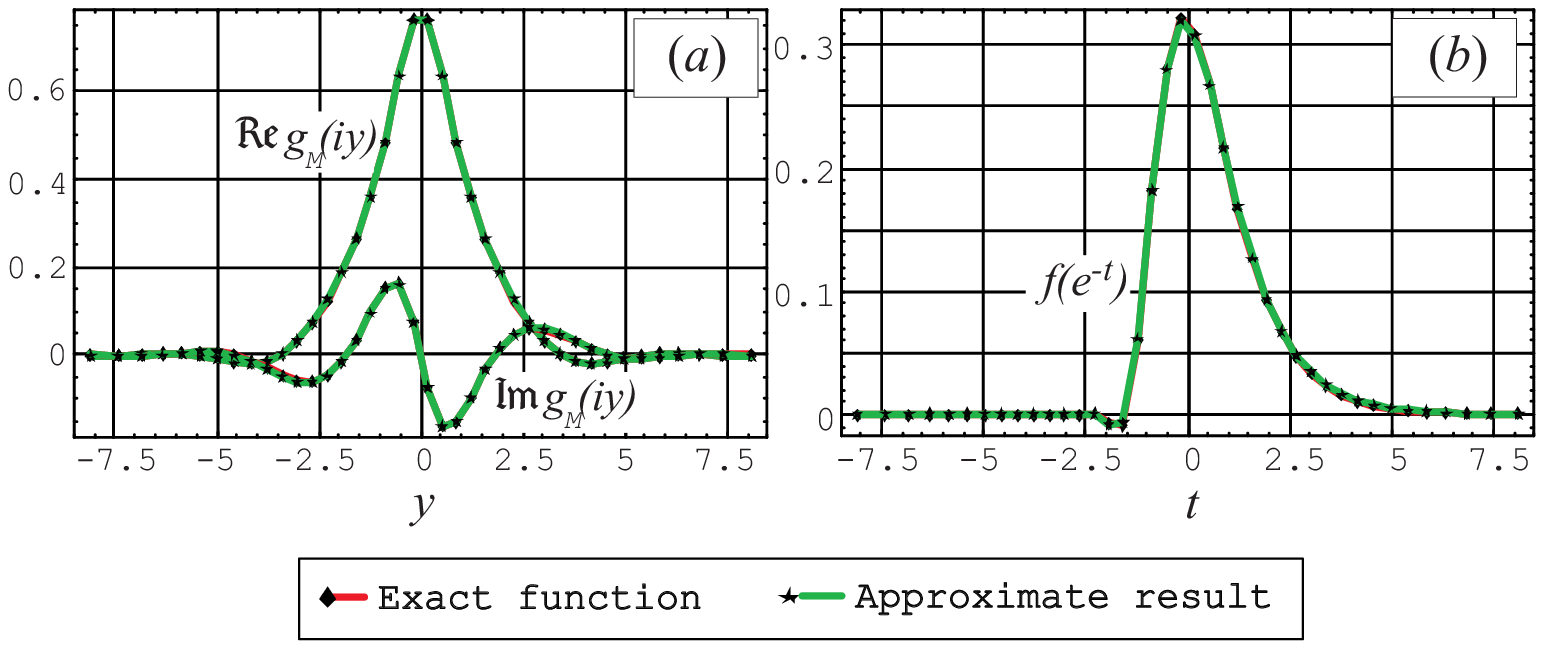}}\hfill}
\begin{center}
\begin{minipage}{13cm}
{\small 
Figure 4: ({\it a}) Exact and approximate Mellin transform of (\ref{ejcmel}) on the imaginary axis. ({\it b}) Inverse transform obtained by (\ref{tdsm}). In both cases 40 Hermite zeros were used.}
\end{minipage}
\end{center}
\section{Final remark}
Finally we address the performance of the discrete Laplace and Mellin transforms presented here on singular problems such as the above first and third examples. Repeated numerical calculations can be done to estimate the convergence of the results yielded by these discrete transforms. Thus, by changing the number $N$ of Hermite zeros it can be seen that in the third example, the Mellin transform of (\ref{ejtmelu}), the relative error goes as $1/\sqrt{N}$. In the case of the first example, the Laplace inversion of (\ref{ejpri}), it is necessary to measure convergence in a different way since it is not possible to evaluate a delta function. To this end, we compute the area under the linear interpolation of the entries of the vector yielded by the numerical Laplace inversion, and test this value against the correct result. For $n=1$ this integral should be 2 and the numerical integrations give 2.0052, 2.0032 and 2.0025, for 50, 80 and 100 Hermite zeros, respectively. In order to give a visual representation of this case, we present in Figure 5 the discrete inverse for $n=1$ and $N=100$.
\vskip1cm
\hbox to \textwidth{\hfill\scalebox{0.8}{\includegraphics{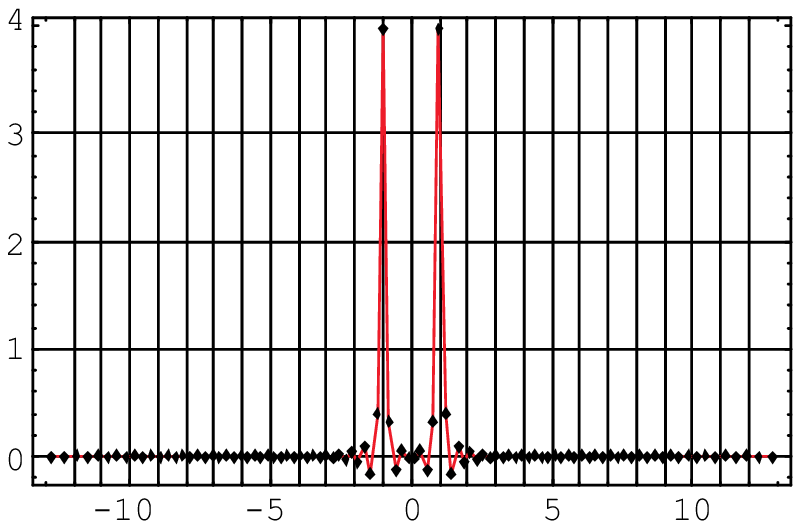}}\hfill}
\begin{center}
\begin{minipage}{13cm}
{\small 
Figure 5: Numerical inversion of the Laplace transformed function (\ref{ejpri}) for $n=1$ and $N=100$.}
\end{minipage}
\end{center}


\end{document}